\documentclass[12pt,reqno]{amsart}

\usepackage{amssymb}

\usepackage{eucal}

\input{diagrams}

\def\BN{{\mathbb N}}

\def\BZ{{\mathbb Z}}

\def\fg{{\mathfrak g}}

\def\fm{{\mathfrak m}}

\def\sD{\mbox{\sf D}}

\def\b{\operatorname{b}}

\def\colim{\operatorname{colim}}
\def\D{\sD}

\def\depth{\operatorname{depth}}

\def\dim{\operatorname{dim}}

\def\Ext{\operatorname{Ext}}

\def\fd{\operatorname{fd}}
\def\fg{\operatorname{fg}}

\def\gr{\mbox{\sf gr}}
\def\Gr{\mbox{\sf Gr}}

\def\h{\operatorname{h}}

\def\inf{\operatorname{inf}}

\def\kfd{\operatorname{k.\!fd\,}}

\def\LTensor{\stackrel{\operatorname{L}}{\otimes}}
\def\opp{\operatorname{op}}
\def\pd{\operatorname{pd}}

\def\RHom{\operatorname{RHom}}

\def\sup{\operatorname{sup}}

\numberwithin{equation}{part}

\newtheorem{Lemma}{Lemma}[section]
\newtheorem{Theorem}[Lemma]{Theorem}

\theoremstyle{definition}
\newtheorem{Definition}[Lemma]{Definition}

\newtheorem{Remark}[Lemma]{Remark}

\newtheorem{Example}[Lemma]{Example}

\def\nc{non-com\-mu\-ta\-ti\-ve}

\def\ABthm{Aus\-lan\-der-Buchs\-baum theorem}

\def\A{A}         

\def\Alm{$\A$-left-mo\-du\-le}

\def\grArm{graded $\A$-right-mo\-du\-le}

\def\DGrA{\D(\A)}
\def\DGrAo{\D(\A^{\opp})}
\def\DbGrA{\D^{\b}(\A)}

\def\DpGrAo{\D^+(\A^{\opp})}

\def\DmfgGrA{\D^-_{\fg}(\A)}

\begin{document}

\title[Ext vanishing]
{Ext vanishing and infinite Auslander-Buchsbaum}

\author{Peter J\o rgensen}
\address{Danish National Library of Science and Medicine, N\o rre
All\'e 49, 2200 K\o \-ben\-havn N, DK--Denmark}
\email{pej@dnlb.dk, www.geocities.com/popjoerg}


\keywords{$k$-flat dimension, depth, infinite \nc\ \ABthm,
Ext groups, vanishing theorem}

\subjclass[2000]{16E30, 16W50}

\begin{abstract} 

A vanishing theorem is proved for Ext groups over \nc\ graded
algebras.  Along the way, an ``infinite'' version is proved
of the \nc\ \ABthm.

\end{abstract}

\maketitle

\setcounter{section}{-1}
\section{Introduction}
\label{sec:introduction}

\noindent
Let $R$ be a noetherian local commutative ring, and let $X$ be a
finitely generated $R$-module of finite projective dimension.  The
classical \ABthm\ states
\[
  \pd X = \depth R - \depth X.
\]
This can also be phrased as an Ext vanishing theorem, namely, if $M$
is any $R$-module, then
\begin{equation}
\label{equ:vanishing}
  \Ext_R^i(X,M) = 0
  \; \mbox{ for } \;
  i > \depth R - \depth X.
\end{equation}

In \cite{Frankild} is proved a surprising variation of this: Suppose
that $R$ is complete in the $\fm$-adic topology.  Then equation
\eqref{equ:vanishing} remains true if $X$ is {\em any} $R$-module of
finite projective dimension, provided $M$ is finitely generated.  In
other words, the condition of being finitely generated is shifted
from $X$ to $M$.

In theorem \ref{thm:vanishing2} below, this result will be
generalized to the situation of a \nc\ noetherian $\BN$-graded
connected algebra.

The route goes through an ``infinite'' version of the \nc\ \ABthm,
given in theorem \ref{thm:AB}.  This result is a substantial
improvement of the original \nc\ \ABthm, as given in
\cite[thm.\ 3.2]{PJIdent}, in that the condition of dealing only with
finitely generated modules is dropped.

The notation of this paper is standard, and is already on record in
several places such as \cite{PJGorHom} or \cite{PJIdent}.  So I will
not say much, except that throughout, $k$ is a field, and $\A$ is a
noetherian $\BN$-graded connected $k$-algebra.  However, let me give
one important word of caution: Everything in sight is graded.  So for
instance, $\DGrA$ stands for $\D(\Gr\,\A)$, the derived category of
the abelian category $\Gr(\A)$ of $\BZ$-graded \Alm s and graded
homomorphisms of degree zero.

\section{Auslander-Buchsbaum}
\label{sec:AB}

\begin{Definition}
For $X$ in $\DGrA$, define the $k$-flat dimension by
\begin{eqnarray*}
  \kfd X & = & -\inf\, k \LTensor_{\A} X.
\end{eqnarray*}
\end{Definition}

\begin{Remark}
\label{rmk:kfdfdpd}
Using a minimal free resolution, it is easy to see that if the
cohomology of $X$ is bounded and finitely generated, then
\[
  \kfd X = \fd X = \pd X,
\]
where $\fd$ stands for flat dimension and $\pd$ stands for projective
dimension. 
\end{Remark}

In the following lemma is used $\DbGrA$, the full subcategory of
$\DGrA$ consisting of complexes with bounded cohomology, and
$\DpGrAo$, the full subcategory of $\DGrAo$ consisting of complexes
with cohomology vanishing in low degrees.

\begin{Lemma}
\label{lem:inf_tensor}
Let $X$ in $\DbGrA$ have $\fd X < \infty$, and let $S$ in
$\DpGrAo$ have $\dim_k \h^i(S) < \infty$ for each $i$.  Then
\[
  \inf\, S \LTensor_{\A} X = \inf\, S - \kfd X.
\]
\end{Lemma}

\begin{proof}
Observe that $\fd X < \infty$ implies 
\[
  \inf\, k \LTensor_{\A} X > -\infty, 
\]
hence 
\begin{equation}
\label{equ:kfd_not_infty}
  \kfd X < \infty.  
\end{equation}

If 
\begin{equation}
\label{equ:infS_infty}
  \inf\, S = \infty
\end{equation}
then $S$ is zero.  Then $S \LTensor_{\A} X$ is also zero, so
\begin{equation}
\label{equ:inf_tensor_infty}
  \inf\, S \LTensor_{\A} X = \infty.
\end{equation}
Combining equations \eqref{equ:kfd_not_infty}, \eqref{equ:infS_infty},
and \eqref{equ:inf_tensor_infty} shows that for $\inf\,S = \infty$,
the lemma's equation reads
\[
  \infty = \infty - \mbox{(something smaller than $\infty$)}
\]
which is true.

So for the rest of the proof, I can assume $\inf\, S < \infty$, that
is, $S$ is non-zero.  Note that $S$ is in $\DpGrAo$, so $\inf\, S >
-\infty$, so $\inf\, S$ is a finite number.

Let me first do the case where $S$ is concentrated in degree zero.
Here $S$ is just a \grArm\ which is finite dimensional over $k$, and
the lemma claims
\begin{equation}
\label{equ:inf_tensor_for_module}
  \inf\, S \LTensor_{\A} X = - \kfd X = \inf\, k \LTensor_{\A} X.
\end{equation}

To prove this, note that $\dim_k S < \infty$ gives that there is a
short exact sequence $0 \rightarrow k(\ell) \longrightarrow S
\longrightarrow \tilde{S} \rightarrow 0$ of \grArm s.  This gives a
distinguished triangle $k(\ell) \longrightarrow S \longrightarrow
\tilde{S} \longrightarrow$ in $\DGrAo$, and tensoring with $X$ and
taking the cohomology long exact sequence gives a sequence consisting
of pieces
\[
  \h^i(k(\ell) \LTensor_{\A} X)
  \longrightarrow \h^i(S \LTensor_{\A} X)
  \longrightarrow \h^i(\tilde{S} \LTensor_{\A} X).
\] 
Induction on $\dim_k S$ now shows that the lowest $i$ with
\[
  \h^i(k \LTensor_{\A} X) \not= 0
\]
equals the lowest $i$ with
\[
  \h^i(S \LTensor_{\A} X) \not= 0,
\]
and equation \eqref{equ:inf_tensor_for_module} follows.  The induction
works because $\inf\, k \LTensor_{\A} X > -\infty$, and even works for
$\inf\, k \LTensor_{\A} X = \infty$.

Let me next do the case where $S$ is general.  There is a standard
spectral sequence 
\[
  E_2^{pq} = \h^p(\h^q(S) \LTensor_{\A} X) 
             \Rightarrow \h^{p+q}(S \LTensor_{\A} X).
\]
Since $\fd X < \infty$, the spectral sequence is first quadrant up to
shift, hence converges.  Now, $\dim_k \h^q(S) < \infty$ holds for each
$q$, so if $\h^q(S)$ is non-zero, then the special case of the lemma
dealt with above applies to $\h^q(S) \LTensor_{\A} X$ and shows
\begin{equation}
\label{equ:inf_equality}
  \inf\,\h^q(S) \LTensor_{\A} X = -\kfd X.
\end{equation}

There are now two cases:  First, if
\begin{equation}
\label{equ:kfd_infty}
  \kfd X = -\infty 
\end{equation}
then equation \eqref{equ:inf_equality} gives that if $\h^q(S)$ is
non-zero, then $\inf\, \h^q(S) \LTensor_{\A} X = \infty$, that is,
$\h^p(\h^q(S) \LTensor_{\A} X)$ is zero for each $p$.  Of course this
also holds for $\h^q(S)$ equal to zero, and so in the spectral
sequence, $E_2^{pq}$ is identically zero.  Therefore the limit
$\h^{p+q}(S \LTensor_{\A} X)$ of the spectral sequence is also zero,
so $S \LTensor_{\A} X$ is zero, so
\begin{equation}
\label{equ:inf_tensor_infty_2}
  \inf\, S \LTensor_{\A} X = \infty.
\end{equation}
But $\inf\, S$ is a finite number, and combining this with equations
\eqref{equ:kfd_infty} and \eqref{equ:inf_tensor_infty_2} says that the
lemma's equation reads
\[
  \infty = \mbox{(a finite number)} - (-\infty)
\]
which is true.

Secondly, if 
\[
  \kfd X > -\infty,
\]
then equation \eqref{equ:inf_equality} says that if $\h^q(S)$ is
non-zero, then $\h^p(\h^q(S) \LTensor_{\A} M)$ is non-zero for $p =
-\kfd X$, but zero for $p < -\kfd X$.  And of course, if $\h^q(S)$ is
zero, then $\h^p(\h^q(S) \LTensor_{\A} X)$ is zero for each $p$.  So
in the spectral sequence, $E_2^{pq}$ is non-zero for $p = -\kfd X$ and
$q = \inf\,S$, but zero for lower $p$ or $q$.  Hence $E_2^{-\kfd
X,\inf\,S}$ can be used in a standard corner argument which shows that
the lowest non-zero term in the limit $\h^{p+q}(S \LTensor_{\A} X)$ of
the spectral sequence has degree $p + q = -\kfd X + \inf\,S$.  Hence
\[
  \inf\, S \LTensor_{\A} X = -\kfd X + \inf\,S,
\]
proving the lemma's equation.
\end{proof}

\begin{Theorem}
[Infinite Auslander-Buchsbaum]
\label{thm:AB}
Assume that $\A$ satisfies $\dim_k \Ext_{\A}^i(k,\A) < \infty$ for each
$i$.  Let $X$ in $\DbGrA$ have $\fd X < \infty$.  Then
\[
  \depth X = \depth \A - \kfd X.
\]
\end{Theorem}

\begin{proof}
I have
\[
  \RHom_{\A}(k,X) \cong \RHom_{\A}(k,\A \LTensor_{\A} X)
  \cong \RHom_{\A}(k,\A) \LTensor_{\A} X,
\]
where the second $\cong$ holds because it holds if $X$ is just a
graded flat module, so also holds if $X$ is quasi-isomorphic to a
bounded complex of graded flat modules, and this is the case since $X$
is in $\DbGrA$ and has $\fd X < \infty$.

Thus
\begin{eqnarray*}
  \depth X & = & \inf\, \RHom_{\A}(k,X) \\
  & = & \inf\, \RHom_{\A}(k,\A) \LTensor_{\A} X \\
  & \stackrel{\rm (a)}{=} & \inf\, \RHom_{\A}(k,\A) - \kfd X \\
  & = & \depth \A - \kfd X,
\end{eqnarray*}
where (a) is by lemma \ref{lem:inf_tensor}.  The lemma applies because
$\RHom_{\A}(k,\A)$ is in $\DpGrAo$, and has $\dim_k \h^i
\RHom_{\A}(k,\A) = \dim_k \Ext_{\A}^i(k,\A) < \infty$ for each $i$ by
assumption.
\end{proof}

\begin{Remark}
\label{rmk:AB}
Theorem \ref{thm:AB} even holds for $\depth \A = \infty$, where the
theorem states $\depth X = \infty$.

On the other hand, suppose $\depth \A < \infty$.  Then it is easy to
see that it makes sense to rearrange the equation in theorem
\ref{thm:AB} as
\[
  \kfd X = \depth \A - \depth X.
\]

If the cohomology of $X$ is bounded and finitely generated, then this
equation reads
\[
  \pd X = \depth \A - \depth X
\]
by remark \ref{rmk:kfdfdpd}.  This is the original \nc\ \ABthm, as
proved in \cite[thm.\ 3.2]{PJIdent}.
\end{Remark}

\section{Ext vanishing}
\label{sec:Ext_vanishing}

\begin{Lemma}
\label{lem:inf_tensor_2}
Let $X$ in $\DbGrA$ have $\fd X < \infty$, and let $T$ in $\DpGrAo$ be
so that $\h^i(T)$ is a graded torsion module for each $i$.  Then
\[
  \inf\, T \LTensor_{\A} X \geq \inf\,T - \kfd X.
\]
\end{Lemma}

\begin{proof}
As in the proof of lemma \ref{lem:inf_tensor}, it is easy to see that
for $\inf\, T = \infty$, the lemma's inequality trivially reads
$\infty \geq \infty$.  So for the rest of the proof I can assume
$\inf\, T < \infty$.  Note that $T$ is in $\DpGrAo$, so $\inf\, T >
-\infty$, so $\inf\, T$ is a finite number.

Let me first do the case where $T$ is concentrated in degree zero.
Here $T$ is just a graded torsion module, and the lemma claims
\[
  \inf\, T \LTensor_{\A} X \geq - \kfd X = \inf\, k \LTensor_{\A} X.
\]

If $F$ is a flat resolution of $X$, then this amounts to
\begin{equation}
\label{equ:to_prove}
  \inf\, T \otimes_{\A} F \geq \inf\, k \otimes_{\A} F.
\end{equation}
To prove this, note that as $T$ is a graded torsion module, it is the
colimit of the system
\[
  T \langle 1 \rangle \subseteq T \langle 2 \rangle \subseteq \cdots
\]
where 
\[
  T \langle j \rangle
  = \{\, t \in T \,\mid\, \A_{\geq j}t = 0 \,\}.
\]
Each quotient $T \langle j \rangle / T \langle j-1 \rangle$ is
annihilated by $\A_{\geq 1}$ so has the form $\coprod_{\alpha}
k(\ell_{\alpha})$, so there are short exact sequences of the form
\[
  0 \rightarrow
  T \langle j-1 \rangle \longrightarrow
  T \langle j \rangle \longrightarrow
  \coprod_{\alpha} k(\ell_{\alpha}) \rightarrow 0.
\]
Tensoring such a sequence with $F$ gives a short exact sequence of
complexes because $F$ consists of graded flat modules.  The corresponding
cohomology long exact sequence consists of pieces
\[
  \h^i(T \langle j-1 \rangle \otimes_{\A} F) \longrightarrow
  \h^i(T \langle j \rangle \otimes_{\A} F) \longrightarrow
  \coprod_{\alpha} \h^i(k \otimes_{\A} F)(\ell_{\alpha}).
\]
Induction on $j$ now makes clear that 
\[
  \h^i(k \otimes_{\A} F) = 0 
\]
implies 
\[
  \h^i(T \langle j \rangle \otimes_{\A} F) = 0 \mbox{ for each $j$,}
\]
and this further gives
\[
  \h^i(T \otimes_{\A} F) \cong
  \h^i(\colim T \langle j \rangle \otimes_{\A} F) \cong
  \colim \h^i(T \langle j \rangle \otimes_{\A} F) = 0,
\]
so the inequality \eqref{equ:to_prove} follows.

Let me next do the case where $T$ is general.  There is a standard
spectral sequence
\[
  E_2^{pq} = \h^p(\h^q(T) \LTensor_{\A} X) 
             \Rightarrow \h^{p+q}(T \LTensor_{\A} X).
\]
Since $\fd X < \infty$, the spectral sequence is first quadrant up to
shift, hence converges.  Now, each $\h^q(T)$ is a graded torsion
module, so if $\h^q(T)$ is non-zero, then the special case of the
lemma dealt with above applies to $\h^q(T) \LTensor_{\A} X$ and shows
\[
  \inf\,\h^q(T) \LTensor_{\A} X \geq -\kfd X.
\]

Hence if $\h^q(T)$ is non-zero, then $\h^p(\h^q(T) \LTensor_{\A} X)$
is zero for $p < -\kfd X$.  And of course, if $\h^q(T)$ is zero, then
$\h^p(\h^q(T) \LTensor_{\A} X)$ is zero for each $p$.  So in the
spectral sequence, $E_2^{pq}$ is zero for $p < -\kfd X$ or $q <
\inf\,T$.  Therefore the limit $\h^{p+q}(T \LTensor_{\A} X)$ of the
spectral sequence is zero for $p + q < -\kfd X + \inf\,T$, and this
says
\[
  \inf\,T \LTensor_{\A} X \geq -\kfd X + \inf\,T,
\]
proving the lemma's inequality.
\end{proof}

\begin{Lemma}
\label{lem:vanishing}
Let $\A$ satisfy $\depth \A < \infty$ and $\dim_k
\Ext_{\A}^i(k,\A) < \infty$ for each $i$.  

Let $X$ in $\DbGrA$ have $\fd X < \infty$, and let $T$ in $\DpGrAo$
be so that $\h^i(T)$ is a graded torsion module for each $i$.  Then
\[
  \inf\, T \LTensor_{\A} X \geq \inf\, T + \depth X - \depth \A.
\]
\end{Lemma}

\begin{proof}
Using lemma \ref{lem:inf_tensor_2} and remark \ref{rmk:AB} gives $\geq$
and $=$ in
\[
  \inf\, T \LTensor_{\A} X \geq
  \inf\, T - \kfd X = \inf\, T + \depth X - \depth \A.
\]
\end{proof}

In the following theorem is used $\DmfgGrA$, the full subcategory
of $\DGrA$ consisting of complexes whose cohomology vanishes in high
degrees and consists of finitely generated graded modules.

\begin{Theorem}
\label{thm:vanishing}
Let $\A$ satisfy $\depth \A < \infty$ and $\dim_k
\Ext_{\A}^i(k,\A) < \infty$ for each $i$.  

Let $X$ in $\DbGrA$ have $\fd X < \infty$, and let $M$ be in
$\DmfgGrA$.  Then
\[
  \sup \RHom_{\A}(X,M) \leq \sup\, M - \depth X + \depth \A.
\]
\end{Theorem}

\begin{proof}
It is easy to see that since $M$ is in $\DmfgGrA$, the Matlis dual
$M^{\prime}$ is in $\DpGrAo$ and has $\h^i(M^{\prime})$ a graded torsion
module for each $i$.  So
\begin{eqnarray*}
  \sup \RHom_{\A}(X,M) 
  & = & \sup \RHom_{\A}(X,M^{\prime \prime}) \\
  & \stackrel{\rm (a)}{=} & \sup ((M^{\prime} \LTensor_{\A} X)^{\prime}) \\
  & = & - \inf\, M^{\prime} \LTensor_{\A} X \\
  & \stackrel{\rm (b)}{\leq} & -\inf\, M^{\prime} - \depth X + \depth \A \\
  & = & \sup\, M - \depth X + \depth \A,
\end{eqnarray*}
where (a) is by adjunction and (b) is by lemma \ref{lem:vanishing}.
\end{proof}

The following is the special case of theorem
\ref{thm:vanishing} where $X$ and $M$ are concentrated in degree zero,
that is, where $X$ and $M$ are graded modules.

\begin{Theorem}
[Ext vanishing]
\label{thm:vanishing2}
Assume that $\A$ satisfies $\depth \A < \infty$ and $\dim_k
\Ext_{\A}^i(k,\A) < \infty$ for each $i$.  

Let $X$ in $\Gr(\A)$ have $\fd X < \infty$, and let $M$ be in
$\gr(\A)$.  Then
\[
  \Ext_{\A}^i(X,M) = 0 
  \; \mbox{ for } \;
  i > \depth \A - \depth X.
\]
\end{Theorem}

This says that for $\fd X < \infty$, the number $\depth \A - \depth
X$ plays the role of projective dimension of $X$, but only with
respect to finitely generated graded modules $M$.  

Of course, this fails when $M$ is general, as illustrated by the
following example.

\begin{Example}
Let $\A$ be the polynomial algebra $k[x]$.  Then the conditions of
theorem \ref{thm:vanishing2} are satisfied, and it is classical that
$\depth \A$ is $1$.

Let $X$ be $k[x,x^{-1}]$.  Then $\depth X \geq 1$ because $X$ is a
graded torsion free module, so $\depth \A - \depth X \leq 0$ and
theorem \ref{thm:vanishing2} gives
\[
  \Ext_{\A}^i(X,M) = 0 
  \; \mbox{ for } \;
  i > 0
\]
for $M$ in $\gr(\A)$.

However, this must fail when $M$ is general, for otherwise $X$ would
be a projective object of $\Gr(\A)$ which it is certainly not.
\end{Example}

\medskip
\noindent
{\bf Acknowledgement. }  
I thank Anders Frankild for showing me \cite{Frankild}.


\begin{thebibliography}{9}









\bibitem{Frankild}  A.\ Frankild, {\it Vanishing of local homology}, to
appear in Math.\ Z.\












\bibitem{PJGorHom}  P.\ J\o rgensen, {\it Gorenstein homomorphisms of
non-commutative rings}, J.\ Algebra {\bf 211} (1999), 240--267.





\bibitem{PJIdent}  P.\ J\o rgensen, {\it Non-commutative graded
homological identities}, J.\ London Math.\ Soc.\ (2) {\bf 57} (1998),
336--350.












\end{thebibliography}
\end{document}